\documentstyle[11pt]{article}

\def\1ox{{ \Omega^1_{\scriptstyle{X}} }}
\def\2ox{{ \Omega^2_{\scriptstyle{X}} }}

\def\ok1{{ \Omega^1_K }}
\def\ok2{{ \Omega^2_K }}
\def\Om{{ \Omega }}
\def\om{{ \omega  }}

\def\O{{ {\cal O} }}

\def\ra{{ \rightarrow }}

\def\a{{ \alpha }}

\def\g{{ \gamma }}

\def\A{{ {\cal A}_{g,n} }}
\def\bA{{ \overline{{\cal A}}_{g,n} }}
\def\deg{{ \mbox{deg} }} 

\def\M{{ {\cal M}_{g,n} }}
\def\bM{{ \overline{{\cal M}}_{g,n} }}

\def\del{{ \nabla }}

\newtheorem{thm}{Theorem}
\newtheorem{cor}{Corollary}
\newtheorem{lem}{Lemma}

\title{ ABC inequalities for some moduli spaces
of log-general type}
\author{Minhyong Kim}
\date{}
\begin{document}
\maketitle

Let $B$ be a smooth projective curve of genus
$\g$ over the complex numbers
and let $f:A \ra B$ be a  non-isotrivial
 semi-abelian scheme over $B$ with
projective generic fiber of relative dimension $g$. 
Let $U\subset B$ be the
locus above which the fibers are projective, and let
$S=B-U$ (a finite set). Thus $f:A_U\ra U$ is abelian,
and $f:A\ra B$ is the connected component of its Neron model.
Denote by $g_0$ the dimension of the fixed part of
$A$ and $s=|S|$.
We will adopt the convention
of using the same notation for the map $f$ and several of
its restrictions, unless an explicit danger of confusion
forces us to do otherwise. Let $e:B\ra A$ be the
identity section, and let $W:=e^*\Om_{A/B}$.
Various authors have dealt with upper and lower bounds for
the degree of $W$.
Faltings \cite{F}, for example, shows that 
$\deg (W) \leq g(3\g+s+1)$
 while Moret-Bailly \cite{MB} shows that
$\deg (W) \leq (g-g_0)(\g-1 )$
in the case where $A/B$ is smooth. Arakelov \cite{A} had earlier
given the bound
$(g-g_0)(\g-1+s/2)$
when $A$ is the connected component of the
Jacobian of a family of stable curves.
In this paper, we improve a bit on Faltings, in the general case: 
\begin{thm} Let $f:A\ra B$ be a non-isotrivial
semi-abelian scheme
of relative dimension $g$ 
with projective generic fiber. Then
$$\deg (W) \leq \frac{(g-g_0)}{2}(2\g-2+s),$$
where $g_0$ is the dimension of the fixed part
and $s$ is the number of non-projective fibers.
\end{thm}
Note that the degree is zero in the isotrivial case, so that
the inequality still holds except for the obvious exception
of $\g=0$.

The method of proof is an easy extension of Moret-Bailly's (and 
Szpiro's \cite{S}, which gives our general result for $g=1$),
and is likely to be known to experts.

The reason it might still be worth writing down in full
is because of the recently emerging connection
 with the ABC conjectures. That is, let $\A$ be the 
moduli space of principally polarized 
abelian varieties of dimension $g$ with full level-$n$
structure. For $n\geq 3$, we have that $\A$
is of log-general type. In fact, according to Mumford (\cite{M}, Proposition
3.4), $\A$ has a toroidal compactification $\A \hookrightarrow \bA$,
such that the compactification divisor $D$ has normal
crossings and if $K$ is the canonical divisor of $\bA$,
then $K+D$ is the pull-back of an ample line bundle on the
Baily-Borel compactification
${\cal A}^*_{g,n}$. The geometric version of the
ABC conjectures deals with maps from curves to
varieties of log-general type \cite{B}.
In our case, if $P:B\ra \bA$ is a map whose image does not
lie in $D$, let $S=P^{-1}(D)$ (this is the inverse image as sets, without
multiplicities) and $U=B-S$. Then there is
a family $f:A\ra B$ of semi-abelian varieties such that
$A_U$ is the abelian family induced by $P|U$.(It also has
a level structure which will be unimportant for
our purposes.) As in the beginning, let $W=e^*\Om_{A/B}$.
Then $(\det W)^{(g+1)}\simeq P^*(K+D)$ (\cite{F}, pp. 339-340).
So as a result, we get
\begin{cor}
$$\deg (P^*(K+D)) \leq \frac{g(g+1)}{2}(2\g-2+s)$$
\end{cor}
This is an inequality of the sort conjectured by
Buium for varieties of log general type, only more precise.
 Vojta has also conjectured
such inequalities at the '98 ABC workshop in Tucson, AZ.

For completeness, we outline how one gets the same kind of
inequalities for the moduli space of curves. That is, let 
$\M$ be the moduli space of curves of genus $g\geq 1$ with
a level $n\geq 3$ structure on its Jacobian. Then $\M$
is of log-general type (\cite{M}, Proposition 4.3). Let
$\bM$ be the compactification of $\M$ constructed by
taking the normalization in $\M$ of the Deligne-Mumford
moduli space of stable curves. Finally, let $Y \ra \bM$ be
a smooth allowable modification (\cite{M}, p. 268),
and $D\subset Y$ be the inverse image of the compactification divisor
of $\bM$. Now, if $P:B\ra Y$ is a map from a curve $B$
whose image does not lie in $D$,  $S=P^{-1}D$, and $U=B-S$,
then there
is a family of stable curves $f:C \ra B$ such that $ C_U$
is the smooth family induced by $P|U$. Let $\om_{C/B}$
be the relative dualizing sheaf.
By Mumford (\cite{M}, proof
of proposition 4.3), we have that $\det (f_*(\om_{C/B}))^{13}
= P^*(K_Y+2D)$. So we get the following corollary of Arakelov's
inequality mentioned above:
$$\deg P^*(K_Y+2D) \leq \frac{13g}{2}(2\g-2+s).$$

\section{Proof of theorem}

Let $f:A \ra B$, $U$, $S$, $s$ be as in the previous section.
We need to prepare the situation a bit. First, let $A_0$
be the fixed-part. Then we may replace $A$ by the connected
component of the Neron model of $A_U/(A_0)_U$, which is
semi-stable since $A$ is, and use the fact that $e^*(\Om_{A_0/B})$
is trivial to reduce to the case where $A$ has no fixed-part
(but the dimension is $g-g_0$).

Next, note that when we make a base-change that is
unramified outside $S$, the quantities $2\g-2+s$ and $\deg (W)$
just multiply by the degree of the base change. (The first
by the Hurwitz formula, and the second by the fact that the
formation of $W=e^*(\Om_{A/B})$ commutes with base change.)
Thus, making  a base-change to the field generated by the three-torsion of
$A$ we may assume
that the $A$ has a level three structure. 

By Chai-Faltings \cite{CF}, theorem VI.1.1,
we have the existence of a proper
smooth variety $X$ with a map $g:X \ra B$ and an
isomorphism $X_U\simeq A_U$, such that $D=g^{-1}(S)$ is
a divisor with normal crossings and with the
crucial property that $$g_*(\Om_{X/B}(\log D)) \simeq
W$$

On the other hand, associated to the family $f:A_U \ra U$ we have the
variation of complex Hodge structures $E_0:=R^1f_*(\Omega^{\cdot}_{X_U/U})$,
which is (among other things) 
a bundle with a connection $$\del_0:E_0\ra E_0\otimes \Om_U.$$
$E_0$ also has the sub-bundle giving the Hodge filtration
$F_0\subset E_0$ which we know to be canonically isomorphic
to $W|U$. Denote by $E$ the Deligne extension of $E_0$
to a bundle on $B$ with a log-connection $$\del:E\ra E\otimes \Om_B(S).$$
Then by Steenbrink \cite{St}, $E$ has a sub-bundle $F$, the saturation
of $F_0$, which can be identified with $g_*(\Om_{X/B}(\log D))$.
Thus, we get an isomorphism $W\simeq F \subset E$. Hence, we
need only bound the degree of $F$.
Denote by $G_0$ the quotient $E_0/F_0$ which therefore extends to
the bundle $G=E/F=R^1g_*(\O_X)$. The (flat) 
polarization $<\cdot, \cdot>$ makes the bundle with connection
 $(E_0,\del_0)$
self-dual, and $F_0$ maximal isotropic, so that there is
also an induced duality between $F_0$ and $G_0$. We will
denote the duality paring between $F_0$ and $G_0$ by the
same brackets as the polarization.
We need the following
\begin{lem}
The bundle $E$ is also self-dual.
\end{lem}
{\em Proof.}
In fact, it is self-dual as a bundle with log-connection. To see this 
recall that $E$ is characterized by
the fact that it has a log-connection $\del$ which extends the
connection on $E_0$ and such that the residue operators on
the points of $S$ have eigenvalues in the interval $[0,1)$.
Since the dual connection $\del^v$ on the dual bundle $E^v$
is defined by $\del^v \phi (e)=d \phi (e)-\phi (\del e)$,
we see that the eigenvalues of the residues of $\del^v$
are the negative of those associated to $\del$. However,
since the abelian variety $A_U$ has semi-stable reduction,
the monodromy operator is unipotent, and thus,
the eigenvalues of the residues are zero. Thus they
are also zero for $\del^v$. But this implies that
$(E^v, \del^v)$ is also a Deligne extension for 
$((E_0)^v, (\del_0)^v)\simeq (E_0,\del_0)$,
so we have an isomorphism $(E,\del) \simeq (E^v,\del^v)$.
\vspace{3mm}

From this lemma and the fact that a saturated subsheaf is
determined by its generic stalk, we get that $F$ and
$G$ are also in duality. 

On $X$ we have the exact sequence of sheaves
$$0\ra g^*\Om_B(S) \ra \Om_X(\log D) \ra \Om_{X/B} (\log D) \ra 0$$
from which we get the log Kodaira-Spencer map
$$\rho:F\ra G\otimes \Om_B(S).$$
As in \cite{MB}, consider the exact sequences
$$ 0\ra \Om_B(S) \ra g_*(\Om_X(\log D)) \ra N \ra 0$$
and
$$ 0 \ra N \ra F \ra G\otimes \Om_B(S).$$
(Where $N$ is defined by these sequences.)

\begin{lem}
 $\deg(N)=0$.
\end{lem}
{\em Proof.}
We will prove this by showing that $N$ is preserved by the
log-connection of $E$. This together with the fact that
the eigenvalues of the residues for the log-connection are
zero, will imply that $\det N$ has a regular connection,
and hence, is of degree zero.

It suffices to check that $N|U \subset F_0 \subset E_0$
is preserved by the connection, or even just to
show that it is generically preserved by the connection.
For this reason,
we will be localizing on the base several times in the argument
to follow, without special mention. For a  section
$\a$ of $\Om_X$, denote by $[a]$ its image in $\Om_{X/B}$.

Recall the usual computation of the Gauss-Manin connection
\cite{K} on the subspace $F_0$ of $E_0$ (We will do the
computation on a section over $U$, the local sheaf theoretic
argument being exactly the same.): Let $\phi \in F_0(U)$ so
that it is an element of $\Om_{X_U/U}(X_U)$. There exists a covering
$\{U_i\}$ of $X_U$ such that $\phi $ is locally liftable to
$\phi_i \in \Om_{X_U}$. Let $t$ be a local parameter on the
base and $v=\frac{d}{dt}$. 
Then $\del_{v }\phi$ has two components, the $F_0$
component is given locally
by the formula $((\del_v \phi )_0)_i=[\a_i]\in \Om_{X_U/U}$, where
$$d\phi_i= \a_i\wedge dt \in 
\Om^2_X$$ (From the fact that $[d\phi]=0$ in
$\Om^2_{X_U/U}$, we get that $d\phi_i$ is of the form $(\cdot)\wedge dt$,
and the $[\a_i]$'s glue together, although
the $\a_i$ may not, in general)
while the $G_0$ component is given by the 1-cocycle
$\{\phi_i (v)-\phi_j (v) \}$. Now, when $\phi \in N$, from the
first exact sequence defining $N$,
there is
a {\em global} (along the fibers)
 lifting $\phi' $ of $\phi$ (after possibly
shrinking $U$), and the $G_0$
component is zero, so we have the global element
$$d\phi' \in \Omega^2_X .$$
 But we also have the contraction
operator $i(v): \Omega^2_X \ra \Om^1_X$, so that we have an element
$i(v)(d\phi') \in \Om_{X_U}$. Locally,  if we  write
$(d\phi')_i=\a_i \wedge dt$, we get that
$i(v)(d\phi')=\a_i(v)dt-\a_i$, and hence, $[-i(v)(d\phi')]=\del_v\phi$,
i.e., $\del_v \phi$ is globally liftable. So $\del_v \phi \in N$,
as was to be shown.

\vspace{3mm}

Now we examine the map $\rho: F \ra G\otimes \Om_B(S)$.

\begin{lem}
$\rho$ factors through the dual of $F/N$:
$$\rho:F \ra (F/N)^v \otimes \Om_B(S) \subset G\otimes \Om_B(S)$$
\end{lem}
{\em Proof.}
Given the duality between $F$ and $G$, the dual of
$F/N$ is exactly the annihilator in $G$ of $N\subset F$.
Since all the subsheaves in question are saturated,
we may check generically that the image of $\rho$
annihilates $N$.
But $\rho$ relates to the Gauss-Manin connection according
to the composition:
$$\rho:F \subset E \stackrel{\del }{\ra}
 E \otimes \Om_B(S) \ra G \otimes \Om_B(S)$$
Now, if $s,t \in F_0$, we have
$<s,t>=0$ so for any local vector field $v$ on $U$,
$$<\del_v s,t>+<s, \del_v t>=v<s,t>=0.$$
Now, the polarization induces the perfect pairing
 between $F_0$ and $G_0$ (by the isotropy of $F_0$). So we have
$$<\rho (s),t>+<s,\rho (t)>=0$$ 
for any $s,t \in F_0$.
Hence, if we assume $s\in N$, then $\rho(s)=0$, so
$<s,\rho(t)>=0$ for all $t\in F_0$,
which is what we want.
\vspace{3mm}

From the lemma, we have an injection
$$F/N \hookrightarrow (F/N)^v\otimes \Om_B(S)$$
of vector bundles of the same rank.
Taking top exterior powers gives us an injection of line bundles,
so
$$\deg(F/N) \leq \deg(F/N)^v+r(2\g-2+s),$$
where $r$ is the rank of $F/N$,
and hence,
$$\deg(F/N) \leq \frac{r}{2}(2\g-2+s)\leq \frac{(g-g_0)}{2}(2\g-2+s)$$
Since $N$ has degree 0, we get
$$\deg( W)=\deg (F) \leq \frac{(g-g_0)}{2}(2\g-2+s)$$
as was to be shown.
\vspace{3mm}

{\bf Acknowledgement:} The author is grateful to Steve Zucker for
discussions on Hodge theory related to this paper, in particular,
for clearing up some confusion in the initial proof of Lemma 1.
This work was supported in part by NSF grant DMS 9701489.

{\footnotesize DEPARTMENT OF MATHEMATICS, UNIVERSITY OF ARIZONA, TUCSON,
AZ 85721, E-MAIL: kim@math.arizona.edu}

\begin{thebibliography}{12}
\bibitem{A}{Arakelov, S.: Families of curves with fixed degeneracy,
Izv. Akad. Nauk. {\bf 35}, 1269-1293 (1971)}

\bibitem{B}{Buium, A.: The ABC theorem for abelian varieties,
Inter. Math. Res. Not. {\bf 5}, 220-233 (1994)}

\bibitem{CF}{Chai, C.-L., Faltings, G.: {\em Degeneration of 
abelian varieties}, Springer Verlag, Berlin-Heidelberg-New York (1990)}
\bibitem{D}{Deligne, P.: Th\'{e}orie de Hodge II, Publ. Math. IHES {\bf 40},
5-58 (1971)}
\bibitem{F}{Faltings, G.: Arakelov's theorem for abelian varieties,
Invent. Math. {\bf 73}, 337-347 (1983)}
\bibitem{K}{Katz, N.: Nilpotent connections and the monodromy theorem;
application of a result of Turrittin, Publ. Math. I.H.E.S. {\bf 39},
355-232 (1970)}
\bibitem{MB}{Moret-Bailly, L.: 
{\em Pinceaux de vari\'{e}t\'{e}s ab\'{e}lian }, Ast\'{e}risque {\bf 129},
Soci\'{e}t\'{e} Math\'{e}matique de France, Paris (1985)}
\bibitem{M}{Mumford, D.: Hirzebruch's proportionality theorem in the
non-compact case, Invent. Math. {\bf 42}, 239-272 (1977)}
\bibitem{St}{Steenbrink, J.: Limits of Hodge structure, Invent. Math. 
{\bf 31}, 229-257 (1976)}
\bibitem{S}{Szpiro, L.: Discriminant et conducteur des courbes
elliptiques,  in {\em S\'{e}minaire sur les pinceaux de courbes elliptiques},
Ast\'{e}risque {\bf 183},
Soci\'{e}t\'{e} Math\'{e}matique de France, Paris (1990) }


\end{thebibliography}
\end{document}